\documentclass{article}
\usepackage{amsmath}
\usepackage{amsthm}
\usepackage{amssymb}
\usepackage{amsfonts}
\usepackage{mathrsfs}
\usepackage{enumerate}
\usepackage{authblk}
\usepackage[normalem]{ulem} 
\usepackage[dvipsnames]{xcolor}
\usepackage[utf8]{inputenc}
\usepackage[english]{babel}
\usepackage[ top = 3 cm, bottom = 2 cm]{geometry} 
\usepackage{graphicx}

\usepackage[colorlinks,citecolor=red,pagebackref,hypertexnames=false]{hyperref}
\usepackage{pdfsync}

\numberwithin{equation}{section}
\theoremstyle{plain}
\newtheorem{theorem}{Theorem}

\newtheorem{corollary}[theorem]{Corollary}

\theoremstyle{definition}

\newtheorem{remark}[theorem]{Remark}
\def\Xint#1{\mathchoice
{\XXint\displaystyle\textstyle{#1}}%
{\XXint\textstyle\scriptstyle{#1}}%
{\XXint\scriptstyle\scriptscriptstyle{#1}}%
{\XXint\scriptscriptstyle\scriptscriptstyle{#1}}%
\!\int}
\def\XXint#1#2#3{{\setbox0=\hbox{$#1{#2#3}{\int}$ }
\vcenter{\hbox{$#2#3$ }}\kern-0.96\wd0}}

\def\dashint{\Xint-}

\begin{document}

\title{Generalized Newton-Leibniz Formula and the Embedding of the Sobolev Functions with Dominating Mixed Smoothness into H\"{o}lder Spaces}
\author{Ugur G. Abdulla\thanks{Ugur.Abdulla@oist.jp} \ }
\affil{Analysis \& PDE Unit, Okinawa Institute of Science and Technology\\ Okinawa, Japan 904-0495}

\maketitle

\abstract{It is well-known that the embedding of the Sobolev space of weakly differentiable functions into H\"{o}lder spaces holds if the integrability exponent is higher than the space dimension. In this paper, the embedding of the Sobolev functions into the H\"{o}lder spaces is expressed in terms of the minimal weak differentiability requirement independent of the integrability exponent. The proof is based on the generalization of the Newton-Leibniz formula to the $n$-dimensional rectangle and inductive application of the Sobolev trace embedding results. The method is applied to prove the embedding of the Sobolev spaces with dominating mixed smoothness into H\"{o}lder spaces.  }



{\large 
\section{Prelude}

Let $W_p^1(\mathbb{R}^n), 1\leq p\leq \infty$ be a Sobolev space of weakly differentiable functions $u\in L_p(Q)$ with first order weak derivatives in $L_p(Q), i=1,...,n$. Originally discovered in the celebrated paper \cite{sobolev}, the concept of Sobolev spaces became a trailblazing idea in many fields of mathematics. The goal of this paper is to analyze embedding of $W_p^1(\mathbb{R}^n)$ into H\"{o}lder spaces $C^{0,\mu}(\mathbb{R}^n), 0\leq \mu \leq 1$ \cite{adams}. Standard notation will be employed for embedding of Banach spaces:
\begin{itemize}
\item $B_1 \hookrightarrow B_2$ means bounded embedding of $B_1$ into $B_2$, i.e. $B_1\subset B_2$, and
\[ \|u\|_{B_2}\leq C \|u\|_{B_1}, \ \forall u\in B_1, \ \text{for some constant} \ C.\]
\item $B_1 \Subset B_2$ denotes compact embedding of $B_1$ into $B_2$, meaning that $B_1 \hookrightarrow B_2$, and every bounded subset of $B_1$ is precompact in $B_2$.
\end{itemize}
If $n=1$, the equivalency class of elements of $W_p^1(\mathbb{R})$ always contain an absolutely continuous element, which is H\"{o}lder continuous with exponent $1-p^{-1}$, if $p>1$,   i.e. there is a bounded embedding
\begin{equation}\label{1dembedding}
W_p^1(\mathbb{R}) \hookrightarrow C^{0,1-\frac{1}{p}}(\mathbb{R}), \ \text{if} \ p>1; \ \ W_1^1(\mathbb{R}) \hookrightarrow C^{0}(\mathbb{R}).
\end{equation}
The embedding \eqref{1dembedding} easily follows from the Newton-Leibniz formula
\begin{equation}\label{newtonleibniz1d}
u(x')-u(x)=\int_x^{x'} \frac{du(y)}{dx}\,dy
\end{equation}
via the application of the H\"{o}lder inequality and compactness argument. The embedding \eqref{1dembedding} fails to be true if $n\geq 2$ and $p\leq n$. However, there is a bounded embedding \cite{morrey}
\begin{equation}\label{ndembedding}
W_p^1(\mathbb{R}^n) \hookrightarrow C^{0,1-\frac{n}{p}}(\mathbb{R}^n) \ \text{if} \ p>n.
\end{equation}
Hence, stretching the integrability exponent $p$ beyond space dimension $n$ implies the H\"{o}lder continuity. In particular, elements of the Hilbert space $H^1(\mathbb{R}^n)=W_2^1(\mathbb{R}^n)$, are not  
continuous in general, if $n\geq 2$. The main goal of this paper is to express the continuity of elements of $W_p^1(\mathbb{R}^n)$ in terms of weak differentiability requirements.

{\bf Problem:} {\it What are the minimal weak differentiability requirements on elements of $W_p^1(\mathbb{R}^n)$ ($1\leq p \leq n$) to be continuous? What is the largest subspace of $W_p^1(\mathbb{R}^n)$ embedded into H\"{o}lder space for all $p>1$?}

The paper reveals that the anticipated subspace is the Sobolev-Nikol'skii space 
\begin{equation*}
    S_p^1(\mathbb{R}^n) = \Big \{ u \in W_p^1(\mathbb{R}^n)\Big | \ \frac{\partial^k u}{\partial x_{i_1} \cdots \partial x_{i_k}}  \in L_p(\mathbb{R}^n), i_1<\cdots <i_k, k=\overline{2,n} \Big \}. 
\end{equation*}
equipped with the norm 
\begin{equation*}
 \|u\|_{S_p^1(\mathbb{R}^n)}:= 
 \left\{
    \begin{array}{l}
\Big ( \|u\|_{L_p(\mathbb{R}^n)}^p+\sum\limits_{k=1}^n\sum\limits_{\substack{i_1,...,i_k=1 \\ i_1<...<i_k}}^n  \Big \| \frac{\partial^k u}{\partial x_{i_1} \cdots \partial x_{i_k}} \Big  \|_{L_p(\mathbb{R}^n)}^p \Big )^{\frac{1}{p}}, \quad\text{if} \ 1\leq p <\infty,\\
\|u\|_{L_\infty(\mathbb{R}^n)}+\sum\limits_{k=1}^n\sum\limits_{\substack{i_1,...,i_k=1 \\ i_1<...<i_k}}^n \Big \| \frac{\partial^k u}{\partial x_{i_1} \cdots \partial x_{i_k}} \Big  \|_{L_\infty(\mathbb{R}^n)}, \quad\text{if} \ p =\infty.
 \end{array}\right.
\end{equation*}
The space $S_p^1(\mathbb{R}^n)$ is a special case of Sobolev spaces with dominating mixed smoothness. The class was introduced by Nikol'skii in \cite{nikolski1,nikolski2}. 
There is a vast literature on analysis of these spaces. We refer to \cite{triebel1,vybiral,triebel2} and the references therein.

The main result of this paper is twofold. First, we introduce and prove a generalization of the celebrated Newton-Leibniz formula to $n$-dimensional rectangles (or $n$-\textit{rectangles}). Then by using the new formula as a tool, we present a surprisingly simple and elegant proof of the embedding of the Sobolev spaces with dominating mixed smoothness into H\"{o}lder spaces. The proof resembles the proof of the embedding \eqref{1dembedding} in the one-dimensional case by using generalized Newton-Leibniz formula, H\"{o}lder inequality and iterative application of the Sobolev trace embedding results.  
\section{Notations}
\label{sec:notations}
\begin{itemize}
\item $C^0(\mathbb{R}^n)$ is a Banach space of continuous and bounded functions with norm
\begin{equation*}
\|u\|_{C^{0}(\mathbb{R}^n)}:=\sup\limits_{x\in\mathbb{R}^n}|u(x)|=\|u\|_{L_\infty(\mathbb{R}^n)}.
\end{equation*}
\item For $k\in \mathbb{N}$, $C^k(\mathbb{R}^n)$ is a Banach space of $k$ times continuously differentiable functions, with all derivatives of order $k$ bounded, and with the norm
\begin{equation*}
\|u\|_{C^{k}(\mathbb{R}^n)}:=\sum\limits_{j=0}^k\sup\limits_{x\in\mathbb{R}^n}|D^ju(x)|=\sum\limits_{j=0}^k\|D^ju\|_{L_\infty(\mathbb{R}^n)},
\end{equation*}
where $D^ju$ is a tensor of rank $j$, dimension $n$, and
\[ |D^ju|=\Big ( \sum\limits_{i_1,...,i_j=1}^n \Big | \frac{\partial^ju(x)}{\partial x_{i_1}\cdots \partial x_{i_j}}\Big |^2\Big )^{\frac{1}{2}} \]
\end{itemize}
The following standard notation will be used for H\"{o}lder spaces:
\begin{itemize}
\item For $0\leq \gamma \leq 1$, H\"{o}lder space $C^{0,\gamma}(\mathbb{R}^n)$ is the Banach space of elements $u\in C^0(\mathbb{R}^n)$ with finite norm
\begin{equation*}
\|u\|_{C^{0,\gamma}(\mathbb{R}^n)}:=\|u\|_{C^0(\mathbb{R}^n)}+[u]_{C^{0,\gamma}(\mathbb{R}^n)}
\end{equation*}
where 
\[ [v]_{C^{0,\gamma}(\mathbb{R}^n)}:= \sup\limits_{\substack{x, x'\in \mathbb{R}^n \\ x\neq x'}} \frac{|v(x)-v(x')|}{|x-x'|^\gamma} \]
The space $C^{0,0}(\mathbb{R}^n)$ is equivalent to $C^{0}(\mathbb{R}^n)$.
\item For $k\in \mathbb{N}$, $0\leq \gamma \leq 1$, H\"{o}lder space $C^{k,\gamma}(\mathbb{R}^n)$ is a subspace of $C^k(\mathbb{R}^n)$ with finite norm
\begin{equation*}
\|u\|_{C^{k,\gamma}(\mathbb{R}^n)}:=\sum\limits_{j=0}^k\|D^ju\|_{C^{0,\gamma}(\mathbb{R}^n)}
\end{equation*}
\end{itemize}
 Throughout the paper we use standard notations for $L_p(Q), 1\leq p \leq \infty$ spaces; the following standard notations are used for Sobolev spaces \cite{adams}:
\begin{itemize}
\item For $k\in \mathbb{N}$, $1\leq p \leq \infty$, Sobolev space $W_p^k((\mathbb{R}^n)$ is the Banach space of measurable functions on $\mathbb{R}^n$ with finite norm
\begin{equation*}
\| u\|_{W_p^k(\mathbb{R}^n)} :=
   \sum\limits_{j=0}^k\|D^ju\|_{L_p(\mathbb{R}^n)}. \end{equation*}
\item For ${\bf s}=(s_1,...,s_n) \in \mathbb{Z}_+^n, 1\leq p \leq \infty$, anisotropic Sobolev space $W_p^{\bf s}(\mathbb{R}^n)$ is the Banach space of measurable functions on $\mathbb{R}^n$ with finite norm
\begin{equation*}
\| u\|_{W_p^{\bf s}(\mathbb{R}^n)} := 
\|u\|_{L_p(\mathbb{R}^n)}+\sum\limits_{i=1}^n\sum\limits_{k=1}^{s_i}  \Big \|\frac{\partial^k u}{\partial x_i^k} \Big \|_{L_p(\mathbb{R}^n)}. 
\end{equation*}
Note that the size of the vector $\bf s$ coincides with the dimension of the space. In particular, for $1\leq k\leq n$, and fixed $j\in \{1,...,k\}$, we consider Sobolev spaces $W_p^{{\bf s}}(\mathbb{R}^k)$ of the weakly $x_j$-differentiable functions on $\mathbb{R}^k$, where ${\bf s}=(s_i)_{i=1}^k \in \mathbb{Z}_+^k$ and $s_i=\delta_{ij}$ is a Kronecker symbol.
\item For ${\bf k}=(k_1,...,k_n)\in \mathbb{Z}_+^n, 1\leq p \leq \infty$, Sobolev space $S_p^{{\bf k}}(\mathbb{R}^n)$ with dominating mixed derivatives is a Banach space of measurable functions on $\mathbb{R}^n$ with finite norm
\begin{equation*}
\| u\|_{S_p^{\bf k}(\mathbb{R}^n)} := 
\sum\limits_{\alpha\in\mathbb{Z}_+^n, \alpha_i\leq k_i}  \Big \|\ \frac{\partial^{|\alpha|}u(x)}{\partial x_1^{\alpha_1}\cdots \partial x_n^{\alpha_n}}\Big \|_{L_p(\mathbb{R}^n)}. 
\end{equation*}
If $k_1=\cdots=k_n=k\in\mathbb{N}$, we shall write $S_p^{k}(\mathbb{R}^n)=S_p^{{\bf k}}(\mathbb{R}^n)$.
\item Let $Q\subset \mathbb{R}^n$ be a bounded domain. For ${\bf k}=(k_1,...,k_n)\in \mathbb{Z}_+^n, 1\leq p \leq \infty$, Sobolev space $S_p^{{\bf k}}(Q)$ with dominating mixed derivatives is defined as 
\[ S_p^{{\bf k}}(Q)=\{f\in {\cal D}'(Q): \ \exists g\in S_p^{{\bf k}}(\mathbb{R}^n) \ \text{with} \ g |_Q=f\} \]
and with
\begin{equation*}
\| u\|_{S_p^{\bf k}(Q)} := 
\inf \| u\|_{S_p^{\bf k}(\mathbb{R}^n)},
\end{equation*}
where the infimum is taken over all $g\in S_p^{\bf k}(\mathbb{R}^n)$ such that its restriction $g|_Q$ to $Q$ coincides with $f$ in the space of distributions ${\cal D}'(Q)$. 

If $k_1=\cdots=k_n=k\in\mathbb{N}$, we shall write $S_p^{k}(Q)=S_p^{{\bf k}}(Q)$.
\end{itemize}
\section{Main Results}
\subsection{Generalized Newton-Leibniz Formula}

Let $x,x'\in \mathbb{R}^n$ with $x_i<x'_i, i=\overline{1,n}$ are fixed and $P$ be $n$-\textit{rectangle}
\begin{equation}\label{nrectangle}
 P=\{\eta\in \mathbb{R}^n: x_i\leq \eta_i\leq x'_i, i=\overline{1,n}\}
 \end{equation}
with vertex $x$ (or $x'$) called a \textit{bottom (or top) corner} of $P$. For any subset $\{i_1,...,i_k\}\subset \{1,...,n\}, k=\overline{1,n}$, let 
\[P_{i_1\dots i_k}=P\cap \{\eta\in \mathbb{R}^n: \eta_l=x_l, l\neq i_j, j=\overline{1,k} \} \]
be a $k$-\textit{rectangle} with bottom corner $x$. Note that $P_{i_1\dots i_k}$ is invariant with respect to permutation of multiindex $i_1\cdots i_k$, and it coincides with $P$ if $k=n$. 

The following is the generalization of the celebrated Newton-Leibniz formula:
\begin{theorem}\label{newtonleibniz}
Any function $u\in C^n(P)$ satisfies the following 
\textbf{generalized Newton-Leibniz formula}:
\begin{equation}\label{NL}
u(x')-u(x)=\sum\limits_{k=1}^n\sum\limits_{\substack{i_1,...,i_k=1 \\ i_1<...<i_k}}^n \ \int\limits_ {P_{i_1\dots i_k}}\frac{\partial^k u(\eta)}{\partial x_{i_1} \cdots \partial x_{i_k}} \,d\eta_{i_1}\cdots\,d\eta_{i_k}.
\end{equation}
\end{theorem}
 If $n=1$, \eqref{NL} coincides with the Newton-Leibniz formula \eqref{newtonleibniz1d}. Note that for $\forall k$ there are ${n}\choose {k}$ integrals in \eqref{NL} along all $k$-\textit{rectangles} $P_{i_1\dots i_k}$ with bottom corner $x$. Therefore, altogether there are 
\[ \sum\limits_{k=1}^n {{n}\choose {k}} = \sum\limits_{k=0}^n {{n}\choose {k}} -1 = (1+1)^n-1=2^n-1 \]
integrals in \eqref{NL} along all sub-rectangles of $P$ with bottom corner at $x$.  
\subsection{Embedding of the Sobolev Spaces with Dominating Mixed Smoothness into H\"{o}lder Spaces}
\begin{theorem}\label{newembedding}
The following bounded embedding holds
\begin{equation}\label{newembd}
S_p^1( \mathbb{R}^n)\hookrightarrow C^{0,1-\frac{1}{p}}( \mathbb{R}^n); \ \  1\leq p \leq \infty,
\end{equation}
The equivalency class of every element of $S_p^1( \mathbb{R}^n)$ possesses a representative in $C^{0,1-\frac{1}{p}}( \mathbb{R}^n)$, which satisfies the 
\textbf{generalized Newton-Leibniz formula} \eqref{NL}, 
where $P\subset  \mathbb{R}^n$ is an $n$-\textit{rectangle} with bottom and top corner at $x$ and $x'$ respectively. In particular, $\forall k=1,...,n-1$ and $1\leq i_1<\cdots <i_k\leq n$
\begin{equation}\label{Lptraces}
\frac{\partial^k u}{\partial x_{i_1} \cdots \partial x_{i_k}} \in L_p(P_{i_1\dots i_k}),
\end{equation}
in the sense of traces.
\end{theorem}
\begin{corollary}
For $k\in \mathbb{N}$ the following bounded embedding holds
\begin{equation}\label{newembd1}
S_p^k( \mathbb{R}^n)\hookrightarrow C^{k-1,1-\frac{1}{p}}( \mathbb{R}^n); \ \  1\leq p \leq \infty,
\end{equation}
\end{corollary}
The following sharp embedding result holds for the anisotropic Sobolev spaces with dominating mixed smoothness:
\begin{corollary}
Let ${\bf k}=(k_1,...,k_n)\in \mathbb{N}^n, 1\leq p \leq \infty$, and $u\in S_p^{{\bf k}}(\mathbb{R}^n)$. 
Then $\forall m=1,...,n$ and $\forall 1\leq i_1<i_2<\cdots <i_m\leq n$
\begin{equation}\label{newembd2}
\frac{\partial^{k_{i_1}+\cdots +k_{i_m}-m}u}{\partial x_{i_1}^{k_{i_1}-1}\cdots \partial x_{i_m}^{k_{i_m}-1}} \in C^{0,1-\frac{1}{p}}( \mathbb{R}^n)
\end{equation}
\end{corollary}
\begin{corollary}
Let $Q\subset \mathbb{R}^n$ be a bounded domain. For $k\in \mathbb{N}$ the following bounded and compact embeddings hold
\begin{equation}\label{newembddom}
S_p^k(Q)\hookrightarrow C^{k-1,1-\frac{1}{p}}(\overline{Q}), \ \text{if} \ 1\leq p \leq \infty;
\end{equation}
\begin{equation}\label{newembdcmp}
S_p^k(Q)\Subset C^{k-1,\mu}(\overline{Q}), 0<\mu < 1-\frac{1}{p}, \ \text{if} \ 1<p\leq \infty;
\end{equation}
\end{corollary}


\section{Proof of Main Results}
\label{proof}
\textit{Proof of Theorem~\ref{newtonleibniz}}. \ Assuming that $u\in C^n(P)$, we prove \eqref{NL} by induction in terms of the space dimension $n$. If $n=1$, it coincides with he Newton-Leibniz formula. Assume that \eqref{NL} is true, and demonstrate that it is true if $n$ is replaced with $n+1$. Let $x,x'\in \mathbb{R}^{n+1}$ with $x_i<x'_i, i=\overline{1,n+1}$, are fixed. We have
\begin{equation}\label{induction1}
u(x')-u(x)=(u(x')-u(\tilde{x},x'_{n+1})) + (u(\tilde{x},x'_{n+1})-u(x)),
\end{equation}
where $\tilde{x}=(x_1,...,x_n)$. Applying \eqref{NL} to the first term and the Newton-Leibniz formula to the second term in \eqref{induction1}, we derive
\begin{gather}\label{induction2}
u(x')-u(x)=\sum\limits_{k=1}^n\sum\limits_{\substack{i_1,...,i_k=1 \\ i_1<...<i_k}}^n \ \int\limits_ {P_{i_1\dots i_k}} \frac{\partial^k u(\tilde{\eta},x'_{n+1})}{\partial x_{i_1} \cdots \partial x_{i_k}} \,d\eta_{i_1}\cdots\,d\eta_{i_k} +
\int\limits_{x_{n+1}}^{x'_{n+1}} \frac{\partial u(\tilde{x},\eta)}{\partial x_{n+1}}\,d\eta.
\end{gather}
Applying Newton-Leibniz formula to all but the last integrand, we have
\begin{gather}
u(x')-u(x)=\sum\limits_{k=1}^n\sum\limits_{\substack{i_1,...,i_k=1 \\ i_1<...<i_k}}^n \ \int\limits_ {P_{i_1\dots i_k}} \int\limits_{x_{n+1}}^{x'_{n+1}}  \frac{\partial^{k+1} u(\eta)}{\partial x_{i_1} \cdots \partial x_{i_k} \partial x_{n+1}} \,d\eta_{i_1}\cdots\,d\eta_{i_k}\,d\eta_{n+1} \nonumber\\ +
\int\limits_{x_{n+1}}^{x'_{n+1}} \frac{\partial u(\tilde{x},\eta)}{\partial x_{n+1}}\,d\eta + \sum\limits_{k=1}^n\sum\limits_{\substack{i_1,...,i_k=1 \\ i_1<...<i_k}}^n \ \int\limits_ {P_{i_1\dots i_k}} \frac{\partial^k u(\eta)}{\partial x_{i_1} \cdots \partial x_{i_k}} \,d\eta_{i_1}\cdots\,d\eta_{i_k},\label{ind3}
\end{gather}
which imply that
\begin{equation}\label{NL-1}
u(x')-u(x)=\sum\limits_{k=1}^{n+1}\sum\limits_{\substack{i_1,...,i_k=1 \\ i_1<...<i_k}}^{n+1} \ \int\limits_ {P_{i_1\dots i_k}} \frac{\partial^k u(\eta)}{\partial x_{i_1} \cdots \partial x_{i_k}} \,d\eta_{i_1}\cdots\,d\eta_{i_k},
\end{equation}
where we use the same notation for the $(n+1)$-\textit{rectangle} $P$, as well as its corresponding sub-rectangles in $\mathbb{R}^{n+1}$. Indeed, divide all $2^{n+1}-1$ sub-rectangles of $P$ with the bottom corner at $x$ into two groups depending on whether or not the edge $p_{n+1}$ joining vertices $x$ and $(\tilde{x},x'_{n+1})$ is contained in it. The first two terms on the right hand side of (\ref{ind3}) consist of all $2^n$ terms of \eqref{NL-1} with integrals along sub-rectangles containing the edge $p_{n+1}$, and the last term on the right hand side of \eqref{ind3} is identical with the remaining $2^n-1$ integrals in \eqref{NL-1} along sub-rectangles which do not contain the edge $p_{n+1}$. This completes the proof by induction.

\textit{Proof of Theorem~\ref{newembedding}}. First, we prove the Theorem assuming that $1\leq p <\infty$. The proof will be pursued in four steps.

{\it Step 1.} Prove that for $u \in S_p^1(\mathbb{R}^n)$, each of the $2^n-1$ integrals on the right hand side of \eqref{NL} is finite, and in particular, \eqref{Lptraces} is satisfied. Existence of the integral with $k=n$ on the right hand side of \eqref{NL} follows from definition of $S_p^1(\mathbb{R}^n)$ and H\"{o}lder inequality. We prove the existence of the remaining $2^n-2$ trace integrals in \eqref{NL} by mathematical induction and  Sobolev trace embedding result. First, we demonstrate that the claim is true if $k=n-1$. Then we show that the claim is true for any $k<n-1$, provided it is true for $k+1$. Indeed, if $k=n-1$, for each of the $n$ integrals we select a unique integer $j$ satisfying
\begin{equation}\label{ejk+1}
j\in \{1,...,n\}\cap \{i_1,...,i_k\}^c
\end{equation}
and define a multi-index ${\bf s}=(s_1,...,s_n) \in \mathbb{Z}_+^n$, where $s_i=\delta_{ij}$ is a Kronecker symbol. We have
\begin{equation}\label{n-1dim}
\frac{\partial^{n-1} u}{\partial x_{i_1} \cdots \partial x_{i_{n-1}}} \in W_p^{\bf s}(P).
\end{equation}
Note that $(n-1)$-\textit{rectangle} $P_{i_1 \cdots i_{n-1}}$ is a boundary of $P$ on the hyperplane $\eta_j=x_j$. Existence of the trace 
\begin{equation}\label{Lptrace}
\frac{\partial^{n-1} u}{\partial x_{i_1} \cdots \partial x_{i_{n-1}}}\in L_p(P_{i_1 \cdots i_{n-1}})
\end{equation}
is a consequence of the Sobolev trace embedding result:
\begin{equation}\label{traceembeddingP}
W_p^{\bf s}(P) \hookrightarrow L_p(P_{i_1 \cdots i_{n-1}}).
\end{equation}
For completeness, we present a proof of \eqref{traceembeddingP}. Consider a function
\begin{equation}\label{zeta1}
\zeta(\eta)=1-\frac{\eta_{j}-x_{j}}{x_{j}'-x_{j}},
\end{equation}
which satisfy
\begin{equation}\label{zetafunction1}
0\leq \zeta \leq 1, \ \ \Big |\frac{\partial \zeta}{\partial \eta_{j}}\Big | \leq \frac{1}{x_j'-x_j}, \ \eta\in P
\end{equation}
Assuming that $u\in C^n(P)$, we have
\begin{gather}
\int\limits_ {P_{i_1\dots i_{n-1}}}\Big | \frac{\partial^{n-1} u(\eta)}{\partial x_{i_1} \cdots \partial x_{i_{n-1}}}\Big |^p \,d\eta_{i_1}\cdots\,d\eta_{i_{n-1}}=\int\limits_ {P_{i_1\dots i_{n-1}}}\zeta\Big | \frac{\partial^{n-1} u(\eta)}{\partial x_{i_1} \cdots \partial x_{i_{n-1}}}\Big |^p \,d\eta_{i_1}\cdots\,d\eta_{i_{n-1}}\nonumber\\
= -\int\limits_ {P_{i_1\dots i_{n-1}}}\int\limits_{x_{j}}^{x_{j}'}\frac{\partial}{\partial x_{j}}\Big ( \zeta \Big | \frac{\partial^{n-1} u(\eta)}{\partial x_{i_1} \cdots \partial x_{i_{n-1}}}\Big |^p\Big ) \,d\eta_{j}\,d\eta_{i_1}\cdots\,d\eta_{i_{n-1}}\nonumber\\
= -\int\limits_ {P_{i_1\dots i_{n-1}}}\int\limits_{x_{j}}^{x_{j}'}\Big [ \frac{\partial\zeta}{\partial x_{j}} \Big | \frac{\partial^{n-1} u(\eta)}{\partial x_{i_1} \cdots \partial x_{i_{n-1}}}\Big |^p+\zeta p \Big | \frac{\partial^{n-1} u(\eta)}{\partial x_{i_1} \cdots \partial x_{i_{n-1}}}\Big |^{p-1} \times \nonumber\\
sgn \Big (  \frac{\partial^{n-1} u(\eta)}{\partial x_{i_1} \cdots \partial x_{i_{n-1}}}\Big ) \frac{\partial^{n} u(\eta)}{\partial x_{i_1} \cdots \partial x_{i_k}\partial x_{j}} \Big ] \,d\eta_{j}\,d\eta_{i_1}\cdots\,d\eta_{i_{n-1}}\label{young1}
\end{gather}
If $p>1$, by using Young's inequality and \eqref{zetafunction1}, from \eqref{young1} it follows 
\begin{equation}\label{Lptraceest}
\Big \|\frac{\partial^{n-1} u(\eta)}{\partial x_{i_1} \cdots \partial x_{i_{n-1}}}\Big \|_{L_p(P_{i_1\cdots i_{n-1}})} \leq C \Big \|\frac{\partial^{n-1} u(\eta)}{\partial x_{i_1} \cdots \partial x_{i_{n-1}}}\Big \|_{W_p^{{\bf s}}(P)},
\end{equation}
where $C=\max(p-1+|x_j'-x_j|^{-1}; 1)$. If $p=1$, \eqref{Lptraceest} follows directly from \eqref{zetafunction1} and \eqref{young1}.
In general, we can approximate $u\in S_p^1(\mathbb{R}^n)$ with the sequence
$u^\epsilon =u*\phi^\epsilon \in C^\infty_{loc}(\mathbb{R}^n)$, where $\phi^\epsilon$ is a standard rescaled mollifier, and derive \eqref{Lptraceest} for $u^\epsilon$. Since $u^\epsilon$ converges to $u$ in the norm given on the right hand side of \eqref{Lptraceest}, it is so in the norm of the left hand side as well, and passing to the limit as $\epsilon \to 0$, \eqref{Lptraceest}, \eqref{traceembeddingP} and \eqref{Lptrace} follow. Hence, $n$ relations of \eqref{Lptraces} with $k=n-1$ are established. Next we prove that the claim is true for $k$, if it is so for $k+1$. For any of the ${{n}\choose {k}}$ integrals in \eqref{NL} along the $k$-dimensional prism $P_{i_1 \cdots i_k}$ we select any integer $j$ satisfying \eqref{ejk+1}, and define a multiindex ${\bf s}=(s_1,...,s_{k+1}) \in \mathbb{Z}_+^{k+1}$, where $s_i=\delta_{ij}$ is a Kronecker symbol. Noting that $P_{i_1 \cdots i_{k+1}}$ is invariant with respect to permutations of the multi-index $i_1 \cdots i_{k+1}$, and due to the induction assumption we have
\begin{equation}\label{kdim}
\frac{\partial^{k} u}{\partial x_{i_1} \cdots \partial x_{i_{k}}} \in W_p^{\bf s}(P_{i_1 \cdots i_kj}).
\end{equation}
$k$-\textit{rectangle} $P_{i_1 \cdots i_{k}}$ is a boundary of $(k+1)$-\textit{rectangle} $P_{i_1 \cdots i_{k} j}$ on the hyperplane $x_j=const$. Sobolev trace embedding result implies:
\begin{equation}\label{traceembeddingP1}
W_p^{\bf s}(P_{i_1 \cdots i_kj}) \hookrightarrow L_p(P_{i_1 \cdots i_{k}}),
\end{equation}
The proof of \eqref{traceembeddingP1} is identical to the proof of \eqref{traceembeddingP}. Hence, \eqref{Lptraces} is proved for all $k$-dimensional integrals. 

{\it Step 2.} In this step we prove that any $u\in S_p^1(\mathbb{R}^n)\cap C^n_{loc}(\mathbb{R}^n), p>1$ satisfies the estimate
\begin{equation}\label{Holder1}
|u(x)-u(x')|\leq \Big [ \Big ((1+p)^{\frac{1}{p}}+ |x-x'|^{\frac{p-1}{p}} \Big )^n - (1+p)^{\frac{n}{p}}\Big ] \|u\|_{S_p^1(\mathbb{R}^n)},
\end{equation}
for all $x,x'\in\mathbb{R}^n$. Similarly, any $u\in S_1^1(\mathbb{R}^n)\cap C^n_{loc}(\mathbb{R}^n)$ satisfy the estimate
\begin{equation}\label{Holder1L1}
|u(x)-u(x')|\leq (3^n-2^n) \|u\|_{S_1^1(\mathbb{R}^n)},
\end{equation}
for all $x,x'\in\mathbb{R}^n$. Note that the estimate \eqref{Holder1L1} is a formal limit of the estimate \eqref{Holder1} as $p\to 1$. 

To prove \eqref{Holder1} (or \eqref{Holder1L1}) without loss of generality we can assume that $x_i<x'_i, i=\overline{1,n}$. Indeed, if $x_i\neq x'_i, i=\overline{1,n}$, then we can 
transform the space via finitely many translations 
\begin{equation}\label{reflection}
\tilde{y}:\mathbb{R}^n\to \mathbb{R}^n, \ \tilde{y}_i=
\left\{
    \begin{array}{l}
y_i, \quad\text{if} \ x_i<x'_i,\\
-y_i, \quad\text{if} \ x_i>x'_i,
 \end{array}\right.
 \end{equation}
and note that the space $S_p^1(\mathbb{R}^n)$ is invariant under this transformation. Then we can apply \eqref{Holder1} (or \eqref{Holder1L1}) to the $\epsilon$-mollification of the transformed function $\tilde{u}(\tilde{x})=u(\tilde{x})$, and passing to limit as $\epsilon \to 0$ deduce \eqref{Holder1} (or \eqref{Holder1L1}) for $\tilde{u}$. Applying inverse transformation of \eqref{reflection} implies \eqref{Holder1} (or \eqref{Holder1L1}) for $u$. If, on the other side $x_i=x'_i$ for some $i$, we can replace $x'_i$ with $x'_i+\delta$, prove \eqref{Holder1} (or \eqref{Holder1L1}) and pass to limit as $\delta\to 0$.

The proof of \eqref{Holder1} and \eqref{Holder1L1} under the assumption that $x_i<x'_i, i=\overline{1,n}$ is based on the generalized Newton-Leibniz formula \eqref{NL}.  
The following is the proof of the estimate \eqref{Holder1}. Let $P$ be a $n$-rectangle \eqref{nrectangle}, and
\[ P^1 :=\{\eta\in\mathbb{R}^n: x_i\leq \eta_i \leq x_i'+1, \ i=\overline{1,n}\} \]

By using H\"{o}lder inequality the integral on the right hand side of \eqref{NL} with $k=n$ is estimated as follows
\begin{equation}\label{est1}
\Big | \int\limits_ {P} \frac{\partial^n u(\eta)}{\partial x_{1} \cdots \partial x_{n}} \,d\eta\Big | \leq |P|^{\frac{p-1}{p}} \Big \|\frac{\partial^n u}{\partial x_{1} \cdots \partial x_{n}}\Big \|_{L_p(P)},
\end{equation}
where $|P|$ denotes volume of the $n$-\textit{rectangle} $P$. For $k=1,...,n-1$, estimation of any of the $k$-dimensional integrals on the right hand side of \eqref{NL} will be pursued in $n-k$ steps. 
Consider typical $k$-dimensional integral in \eqref{NL} along the $k$-\textit{rectangle} $P_{i_1 \cdots i_{k}}$. The idea is based on successive application of the trace embedding result \eqref{traceembeddingP1} $n-k$ times. First we select any integer $j$ from \eqref{ejk+1}, and assign it to multiindex component $i_{k+1}$. By using H\"{o}lder inequality we have
\begin{equation}\label{est2}
\Big | \int\limits_{P_{i_1 \cdots i_k}} \frac{\partial^k u(\eta)}{\partial x_{i_1} \cdots \partial x_{i_k}} \,d\eta_{i_1}\cdots\,d\eta_{i_k}\Big | \leq |P_{i_1 \cdots i_k}|^{\frac{p-1}{p}} \Big \|\frac{\partial^k u}{\partial x_{i_1} \cdots \partial x_{i_k}}\Big \|_{L_p(P_{i_1 \cdots i_k})}
\end{equation}
Consider a function
\begin{equation}\label{zeta}
\zeta(\eta)=1-\frac{\eta_{i_{k+1}}-x_{i_{k+1}}}{x_{i_{k+1}}'-x_{i_{k+1}}+1},
\end{equation}
which satisfy
\begin{equation}\label{zetafunction}
0\leq \zeta \leq 1, \ \ \Big |\frac{\partial \zeta}{\partial \eta_{i_{k+1}}}\Big | \leq 1.
\end{equation}
We have
\begin{gather}
\int\limits_ {P_{i_1\dots i_k}}\Big | \frac{\partial^k u(\eta)}{\partial x_{i_1} \cdots \partial x_{i_k}}\Big |^p \,d\eta_{i_1}\cdots\,d\eta_{i_k}=\int\limits_ {P_{i_1\dots i_k}}\zeta\Big | \frac{\partial^k u(\eta)}{\partial x_{i_1} \cdots \partial x_{i_k}}\Big |^p \,d\eta_{i_1}\cdots\,d\eta_{i_k}\nonumber\\
= -\int\limits_ {P_{i_1\dots i_k}}\int\limits_{x_{i_{k+1}}}^{x_{i_{k+1}}'+1}\frac{\partial}{\partial x_{i_{k+1}}}\Big ( \zeta \Big | \frac{\partial^k u(\eta)}{\partial x_{i_1} \cdots \partial x_{i_k}}\Big |^p\Big ) \,d\eta_{i_{k+1}}\,d\eta_{i_1}\cdots\,d\eta_{i_k}\nonumber\\
= -\int\limits_ {P_{i_1\dots i_k}}\int\limits_{x_{i_{k+1}}}^{x_{i_{k+1}}'+1}\Big [ \frac{\partial\zeta}{\partial x_{i_{k+1}}} \Big | \frac{\partial^k u(\eta)}{\partial x_{i_1} \cdots \partial x_{i_k}}\Big |^p+\zeta p \Big | \frac{\partial^k u(\eta)}{\partial x_{i_1} \cdots \partial x_{i_k}}\Big |^{p-1} \times \nonumber\\
sgn \Big (  \frac{\partial^k u(\eta)}{\partial x_{i_1} \cdots \partial x_{i_k}}\Big ) \frac{\partial^{k+1} u(\eta)}{\partial x_{i_1} \cdots \partial x_{i_k}\partial x_{i_{k+1}}} \Big ] \,d\eta_{i_{k+1}}\,d\eta_{i_1}\cdots\,d\eta_{i_k}\label{young}
\end{gather}
By using Young's inequality and \eqref{zetafunction}, from \eqref{young} it follows 
\begin{gather}
\int\limits_ {P_{i_1\cdots i_k}}\Big | \frac{\partial^k u(\eta)}{\partial x_{i_1} \cdots \partial x_{i_k}}\Big |^p \,d\eta_{i_1}\cdots\,d\eta_{i_k} \leq \int\limits_ {P_{i_1\cdots i_k}}\int\limits_{x_{i_{k+1}}}^{x_{i_{k+1}}'+1} \Big [  \Big |\frac{\partial^{k+1} u(\eta)}{\partial x_{i_1} \cdots \partial x_{i_k}\partial x_{i_{k+1}}}\Big |^p \nonumber\\
 +p \ \Big | \frac{\partial^k u(\eta)}{\partial x_{i_1} \cdots \partial x_{i_k}}\Big |^p\Big ]\,d\eta_{i_{k+1}}\,d\eta_{i_1}\cdots \,d\eta_{i_k}. \label{k+1dimint}
\end{gather}
From \eqref{est2},\eqref{k+1dimint} it follows that
\begin{gather}
\Big | \int\limits_{P_{i_1 \cdots i_k}} \frac{\partial^k u(\eta)}{\partial x_{i_1} \cdots \partial x_{i_k}} \,d\eta_{i_1}\cdots\,d\eta_{i_k}\Big | \leq |P_{i_1 \cdots i_k}|^{\frac{p-1}{p}} \times \nonumber\\
\Big ( \Big \|\frac{\partial^{k+1} u}{\partial x_{i_1} \cdots \partial x_{i_{k+1}}}\Big \|^p_{L_p(P^1_{i_1 \cdots i_k})} + p \Big \|\frac{\partial^k u}{\partial x_{i_1} \cdots \partial x_{i_k}}\Big \|^p_{L_p(P^1_{i_1 \cdots i_k})} \Big )^{\frac{1}{p}}\label{step1est}
\end{gather}
where $P^1_{i_1 \cdots i_k}=P_{i_1 \cdots i_k}\times (x_{i_{k+1}},x_{i_{k+1}}'+1)$ is a $(k+1)$-rectangle. This completes one out of $n-k$ steps for the estimation of the $k$-dimensional integral in \eqref{NL} along the $k$-\textit{rectangle} $P_{i_1 \cdots i_{k}}$. In the next step we select any integer $j$ from \eqref{ejk+1} with $k$ replaced with $k+1$, and assign it to multiindex component $i_{k+2}$. Then for each of the $k+1$-dimensional integrals on the right-hand side of \eqref{step1est} we derive the estimation similar to \eqref{k+1dimint}, where 
$P_{i_1 \cdots i_k}$ is replaced with $(k+1)$-\textit{rectangle} $P^1_{i_1 \cdots i_k}$, and integration interval $(x_{i_{k+1}},x_{i_{k+1}}'+1)$ is replaced accordingly with  $(x_{i_{k+2}},x_{i_{k+2}}'+1)$. Application of these estimations to the right hand side of \eqref{step1est} would complete the second out of $n-k$ steps. By repeating the procedure after $m=1,...,n-k$ steps we derive the following estimate:
\begin{gather}
\Big | \int\limits_{P_{i_1 \cdots i_k}} \frac{\partial^k u(\eta)}{\partial x_{i_1} \cdots \partial x_{i_k}} \,d\eta_{i_1}\cdots\,d\eta_{i_k}\Big | \leq |P_{i_1 \cdots i_k}|^{\frac{p-1}{p}} \times \nonumber\\
\Big [ \sum\limits_{j=0}^{m} {{m}\choose {j}} p^j  \Big \|\frac{\partial^{k+m-j} u}{\partial x_{i_1} \cdots \partial x_{i_{k+m-j}}}\Big \|^p_{L_p(P^{m}_{i_1 \cdots i_k})} \Big ]^{\frac{1}{p}}, \label{stepn-kest}
\end{gather}
where
\[P^m_{i_1 \cdots i_k}=P_{i_1 \cdots i_k}\times (x_{i_{k+1}},x_{i_{k+1}}'+1)\times \cdots \times (x_{i_{k+m}},x_{i_{k+m}}'+1)\]
be a $(k+m)$-\textit{rectangle}. Let us prove the estimation \eqref{stepn-kest} by induction. If $m=1$, the estimation \eqref{stepn-kest} coincides with \eqref{step1est}. Prove that \eqref{stepn-kest} is true for $m+1$ if it is so for any $m<n-k$. Each of the $k+m$-dimensional integrals on the right-hand side of \eqref{stepn-kest} satisfy the following estimate
\begin{gather}
\int\limits_ {P^m_{i_1\cdots i_{k}}}\Big | \frac{\partial^{k+m-j} u(\eta)}{\partial x_{i_1} \cdots \partial x_{i_{k+m-j}}}\Big |^p \,d\eta_{i_1}\cdots\,d\eta_{i_{k+m}} \leq\nonumber\\ \int\limits_ {P^m_{i_1\cdots i_{k}}}\int\limits_{x_{i_{k+1+m}}}^{x_{i_{k+1+m}}'+1} \Big [  \Big |\frac{\partial^{k+1+m-j} u(\eta)}{\partial x_{i_1} \cdots \partial x_{i_{k+m-j}}\partial x_{i_{k+1+m-j}}}\Big |^p +\nonumber\\
 p \ \Big | \frac{\partial^{k+m-j} u(\eta)}{\partial x_{i_1} \cdots \partial x_{i_{k+m-j}}}\Big |^p\Big ]\,d\eta_{i_{k+1+m}}\,d\eta_{i_1}\cdots \,d\eta_{i_{k+m}}. \label{k+mdimint}
\end{gather}
Using \eqref{k+mdimint}, we have
\begin{gather}
\sum\limits_{j=0}^{m} {{m}\choose {j}} p^j  \Big \|\frac{\partial^{k+m-j} u}{\partial x_{i_1} \cdots \partial x_{i_{k+m-j}}}\Big \|^p_{L_p(P^{m}_{i_1 \cdots i_k})}\leq \nonumber\\
\sum\limits_{j=0}^{m} {{m}\choose {j}} p^j  \Big [ \Big \|\frac{\partial^{k+1+m-j} u}{\partial x_{i_1} \cdots \partial x_{i_{k+1+m-j}}}\Big \|^p_{L_p(P^{m+1}_{i_1 \cdots i_k})} + p \Big \|\frac{\partial^{k+m-j} u}{\partial x_{i_1} \cdots \partial x_{i_{k+m-j}}}\Big \|^p_{L_p(P^{m+1}_{i_1 \cdots i_k})} \Big ] =\nonumber\\
\sum\limits_{j=0}^{m} {{m}\choose {j}} p^j  \Big \|\frac{\partial^{k+1+m-j} u}{\partial x_{i_1} \cdots \partial x_{i_{k+1+m-j}}}\Big \|^p_{L_p(P^{m+1}_{i_1 \cdots i_k})}+\nonumber\\\sum\limits_{j=1}^{m+1} {{m}\choose {j-1}} p^j  \Big \|\frac{\partial^{k+1+m-j} u}{\partial x_{i_1} \cdots \partial x_{i_{k+1+m-j}}}\Big \|^p_{L_p(P^{m+1}_{i_1 \cdots i_k})}. \label{m>m+1ind}
\end{gather}
Since
\[ {{m}\choose {j}} +{{m}\choose {j-1}} = {{m+1}\choose {j}}, \ j=1,...,m \]
from \eqref{m>m+1ind} it follows 
\begin{gather}
\sum\limits_{j=0}^{m} {{m}\choose {j}} p^j  \Big \|\frac{\partial^{k+m-j} u}{\partial x_{i_1} \cdots \partial x_{i_{k+m-j}}}\Big \|^p_{L_p(P^{m}_{i_1 \cdots i_k})}\leq \nonumber\\
\sum\limits_{j=0}^{m+1} {{m+1}\choose {j}} p^j  \Big \|\frac{\partial^{k+m+1-j} u}{\partial x_{i_1} \cdots \partial x_{i_{k+1+m-j}}}\Big \|^p_{L_p(P^{m+1}_{i_1 \cdots i_k})},
\end{gather}
which completes the proof of \eqref{stepn-kest} by mathematical induction. 
By choosing $m=n-k$ in \eqref{stepn-kest}, we derive an upper bound of the right hand side of by replacing integration domain with $\mathbb{R}^n$:
 \begin{gather}
\Big | \int\limits_{P_{i_1 \cdots i_k}} \frac{\partial^k u(\eta)}{\partial x_{i_1} \cdots \partial x_{i_k}} \,d\eta_{i_1}\cdots\,d\eta_{i_k}\Big | \leq |P_{i_1 \cdots i_k}|^{\frac{p-1}{p}} \Big [ \sum\limits_{j=0}^{n-k} {{n-k}\choose {j}} p^j \Big ]^{\frac{1}{p}} \|u\|_{S_p^1(\mathbb{R}^n)}\nonumber\\
\leq |x-x'|^{\frac{k(p-1)}{p}} (1+p)^{\frac{n-k}{p}}  \|u\|_{S_p^1(\mathbb{R}^n)}.\label{est3}
\end{gather}
Note that the estimation \eqref{est3} holds for k=n as well in view of \eqref{est1}. By using \eqref{est3} from the generalized Newton-Leibniz formula \eqref{NL} it follows
the estimate
\begin{gather}
|u(x')-u(x)|\leq \sum\limits_{k=1}^n\sum\limits_{\substack{i_1,...,i_k=1 \\ i_1<...<i_k}}^n \ |x-x'|^{\frac{k(p-1)}{p}} (1+p)^{\frac{n-k}{p}}  \|u\|_{S_p^1(\mathbb{R}^n)} \nonumber\\
= \sum\limits_{k=1}^n  \  {{n}\choose {k}} |x-x'|^{\frac{k(p-1)}{p}} (1+p)^{\frac{n-k}{p}}  \|u\|_{S_p^1(\mathbb{R}^n)}\nonumber\\
 = \Big [ \sum\limits_{k=0}^n  \  {{n}\choose {k}} |x-x'|^{\frac{k(p-1)}{p}} (1+p)^{\frac{n-k}{p}} - (1+p)^{\frac{n}{p}}\Big ] \|u\|_{S_p^1(\mathbb{R}^n)}\nonumber\\
 = \Big [ \Big ((1+p)^{\frac{1}{p}}+ |x-x'|^{\frac{p-1}{p}} \Big )^n - (1+p)^{\frac{n}{p}}\Big ] \|u\|_{S_p^1(\mathbb{R}^n)},\label{EST}
\end{gather}
which proves the desired estimate \eqref{Holder1}. The proof of the estimate \eqref{Holder1L1} is almost identical to the proof of \eqref{Holder1}.

{\it Step 3.} In this step we prove 
\begin{itemize}
\item the uniform $C^{0,1-\frac{1}{p}}( \mathbb{R}^n)$-estimate for any $u\in S_p^1(\mathbb{R}^n)\cap C^n_{loc}(\mathbb{R}^n), 1<p<\infty$;
\item the uniform $C^{0}( \mathbb{R}^n)$-estimate for any $u\in S_1^1(\mathbb{R}^n)\cap C^n_{loc}(\mathbb{R}^n)$;
\end{itemize}
Assume $p>1$ and fix $x, x'\in \mathbb{R}^n$ such that $|x-x'|\leq 1$. From \eqref{Holder1} it follows that
\begin{gather}
|u(x)-u(x')|\leq \Big [ \Big ((1+p)^{\frac{1}{p}}+ |x-x'|^{\frac{p-1}{p}} \Big )^n - (1+p)^{\frac{n}{p}}\Big ] \|u\|_{S_p^1(\mathbb{R}^n)}=\nonumber\\
\sum\limits_{k=1}^n  \  {{n}\choose {k}} |x-x'|^{\frac{k(p-1)}{p}} (1+p)^{\frac{n-k}{p}}  \|u\|_{S_p^1(\mathbb{R}^n)}=\nonumber\\
\sum\limits_{k=1}^n  \  {{n}\choose {k}} |x-x'|^{\frac{(k-1)(p-1)}{p}} (1+p)^{\frac{n-k}{p}}  \|u\|_{S_p^1(\mathbb{R}^n)}|x-x'|^{\frac{p}{p-1}}\leq \nonumber\\
\sum\limits_{k=1}^n  \  {{n}\choose {k}} (1+p)^{\frac{n-k}{p}}  \|u\|_{S_p^1(\mathbb{R}^n)}|x-x'|^{\frac{p}{p-1}}=\nonumber\\ 
\Big [\sum\limits_{k=0}^n  \  {{n}\choose {k}} (1+p)^{\frac{n-k}{p}}- (1+p)^{\frac{n}{p}}\Big ] \|u\|_{S_p^1(\mathbb{R}^n)}|x-x'|^{\frac{p}{p-1}}=\nonumber\\
\Big [ \Big (1+(1+p)^{\frac{1}{p}}\Big )^n-(1+p)^{\frac{n}{p}}\Big]  \|u\|_{S_p^1(\mathbb{R}^n)}|x-x'|^{\frac{p}{p-1}}.\label{locHolder}
\end{gather}
Hence, we have
\begin{equation}\label{locHolder1}
\sup\limits_{\substack{|x-x'|\leq 1 \\ x\neq x'}} \frac{|u(x)-u(x')|}{|x-x'|^{\frac{p}{p-1}}}\leq \Big [ \Big (1+(1+p)^{\frac{1}{p}}\Big )^n-(1+p)^{\frac{n}{p}}\Big]  \|u\|_{S_p^1(\mathbb{R}^n)}.
\end{equation}
Now fix $x\in\mathbb{R}^n$. By using \eqref{locHolder1} and 
H\"{o}lder inequality we deduce
\begin{gather}
|u(x)|\leq \dashint\limits_{|y-x|\leq 1}|u(x)-u(y)|\,dy + \dashint\limits_{|y-x|\leq 1}|u(y)|\,dy \leq \nonumber\\
 \Big [ \Big (1+(1+p)^{\frac{1}{p}}\Big )^n-(1+p)^{\frac{n}{p}}\Big]  \|u\|_{S_p^1(\mathbb{R}^n)}+\Gamma_n^{-\frac{1}{p}}\|u\|_{L_p(\mathbb{R}^n)} \leq \nonumber\\
  \Big [ \Big (1+(1+p)^{\frac{1}{p}}\Big )^n-(1+p)^{\frac{n}{p}}+\Gamma_n^{-\frac{1}{p}}\Big]  \|u\|_{S_p^1(\mathbb{R}^n)}.\label{maxnorm}
\end{gather}
where $\Gamma_n$ is a volume of the unit ball in $\mathbb{R}^n$. Hence, we have
\begin{equation}\label{C0normbound}
\|u\|_{C^{0}(\mathbb{R}^n)}\leq \Big [ \Big (1+(1+p)^{\frac{1}{p}}\Big )^n-(1+p)^{\frac{n}{p}}+\Gamma_n^{-\frac{1}{p}}\Big]  \|u\|_{S_p^1(\mathbb{R}^n)}.
\end{equation}
From \eqref{C0normbound} it follows that
\begin{equation}\label{locHolder2}
\sup\limits_{\substack{|x-x'|\geq 1}} \frac{|u(x)-u(x')|}{|x-x'|^{\frac{p}{p-1}}}\leq 2 \|u\|_{C^{0}(\mathbb{R}^n)} \leq 2\Big [ \Big (1+(1+p)^{\frac{1}{p}}\Big )^n-(1+p)^{\frac{n}{p}}+\Gamma_n^{-\frac{1}{p}}\Big]  \|u\|_{S_p^1(\mathbb{R}^n)}.
\end{equation}
From \eqref{locHolder1} and \eqref{locHolder2} we deduce the following H\"{o}lder seminorm estimate for $u$:
\begin{equation}\label{seminormest}
[u]_{C^{0,1-\frac{1}{p}}(\mathbb{R}^n)} \leq 2\Big [ \Big (1+(1+p)^{\frac{1}{p}}\Big )^n-(1+p)^{\frac{n}{p}}+\Gamma_n^{-\frac{1}{p}}\Big]  \|u\|_{S_p^1(\mathbb{R}^n)}.
\end{equation}
Finally, \eqref{C0normbound}, \eqref{seminormest} imply the following H\"{o}lder norm estimate for $u\in S_p^1(\mathbb{R}^n)\cap C^n_{loc}(\mathbb{R}^n), 1<p<\infty$:
\begin{equation}\label{Holdernormest}
\|u\|_{C^{0,1-\frac{1}{p}}(\mathbb{R}^n)} \leq 3\Big [ \Big (1+(1+p)^{\frac{1}{p}}\Big )^n-(1+p)^{\frac{n}{p}}+\Gamma_n^{-\frac{1}{p}}\Big]  \|u\|_{S_p^1(\mathbb{R}^n)}.
\end{equation}
If $p=1$ from the estimate \eqref{Holder1L1} with the similar argument as in \eqref{maxnorm} we derive the following $C^{0}(\mathbb{R}^n)$-estimate for any $u\in S_1^1(\mathbb{R}^n)\cap C^n_{loc}(\mathbb{R}^n)$:
\begin{equation}\label{HoldernormestL1}
\|u\|_{C^{0}(\mathbb{R}^n)} \leq \Big [ 3^n-2^n+\Gamma_n^{-1}\Big]  \|u\|_{S_1^1(\mathbb{R}^n)}.
\end{equation}

{\it Step 4.} We complete the proof of the embedding \eqref{newembd} by using estimates \eqref{Holdernormest}, \eqref{HoldernormestL1} and smooth approximation of elements of $S_p^1(\mathbb{R}^n)$. Given $u\in S_p^1(\mathbb{R}^n), 1\leq p<\infty$, we select a sequence $v_m\in C_0^\infty(\mathbb{R}^n)$ such that
\begin{equation}\label{approximation}
\|v_m- u\|_{S_p^1(\mathbb{R}^n)} \to 0, \ \text{as} \ m \to \infty.
\end{equation}
For example, the sequence $v_m$ can be given explicitly as in \cite{tao} (Lemma 23):
\[ v_m(x)=u^{\frac{1}{m}}(x) \eta \Big (\frac{x}{m}\Big ), \]
where $u^{\frac{1}{m}} =u*\phi^{\frac{1}{m}} \in C^\infty_{loc}(\mathbb{R}^n)\cap S_p^1(\mathbb{R}^n)$ is the $\frac{1}{m}$-mollification of $u$, $\phi^{\frac{1}{m}}$ is a standard rescaled mollifier, 
$\eta \in C_0^\infty(\mathbb{R}^n)$ be a compactly supported function which equals 1 near the origin. 
If $p>1$, then by applying the estimate \eqref{Holdernormest} to $v_m$, we have
\begin{equation}\label{holderestun}
\|v_m\|_{C^{0,1-\frac{1}{p}}(\mathbb{R}^n)} \leq 3\Big [ \Big (1+(1+p)^{\frac{1}{p}}\Big )^n-(1+p)^{\frac{n}{p}}+\Gamma_n^{-\frac{1}{p}}\Big]  \|v_m\|_{S_p^1(\mathbb{R}^n)}.
\end{equation}
Equivalently, we have
\begin{equation}\label{cauchyseq}
\|v_m-v_l\|_{C^{0,1-\frac{1}{p}}(\mathbb{R}^n)} \leq 3\Big [ \Big (1+(1+p)^{\frac{1}{p}}\Big )^n-(1+p)^{\frac{n}{p}}+\Gamma_n^{-\frac{1}{p}}\Big]  \|v_m-v_l\|_{S_p^1(\mathbb{R}^n)}.
\end{equation}
for all $m,l\geq 1$, whence there exists a function $u_*\in C^{0,1-\frac{1}{p}}(\mathbb{R}^n)$ such that
\begin{equation}\label{Holderconv}
\|v_m-u_*\|_{C^{0,1-\frac{1}{p}}(\mathbb{R}^n)} \to 0, \ \text{as} \ m \to \infty.
\end{equation}
From \eqref{approximation} it follows that $u_*=u$, a.e. on $\mathbb{R}^n$, so that $u_*$ is in the equivalency class of $u$. Passing to limit as $m\to \infty$, from \eqref{holderestun} it also follows that
\begin{equation}\label{holderestu*}
\|u_*\|_{C^{0,1-\frac{1}{p}}(\mathbb{R}^n)}\leq 3\Big [ \Big (1+(1+p)^{\frac{1}{p}}\Big )^n-(1+p)^{\frac{n}{p}}+\Gamma_n^{-\frac{1}{p}}\Big]  \|u\|_{S_p^1(\mathbb{R}^n)}.
\end{equation}
which proves the bounded embedding \eqref{newembd}. Step 1 of the proof implies that the traces of $u_*$ satisfy \eqref{Lptraces}, and each of them is an $L_p$-limit of the corresponding sequence of traces of $v_m$. Therefore, writing \eqref{NL} for $v_m$, and passing to limit as $m\to\infty$, it follows that $u_*$ satisfies the generalized Newton-Leibniz formula \eqref{NL}.   

Proof in the case $p=1$ is identical by using an estimate \eqref{HoldernormestL1}. This completes the proof of the theorem in the case $1\leq p <\infty$. 

Assume that $p=\infty$. In this case the embedding \eqref{newembd} is not new, and it is contained in the well-known fact that \cite{adams}
\[ W_\infty^1(\mathbb{R}^n)  \hookrightarrow C^{0,1}(\mathbb{R}^n). \]
Since $S_\infty^1(\mathbb{R}^n)$ is a subspace of $W_\infty^1(\mathbb{R}^n)$, its elements are bounded and Lipschitz continuous functions, and the embedding \eqref{newembd} holds. The assertion that $u$ satisfies \eqref{NL} follows from the proof given for the case $p<\infty$. It only remains to show that \eqref{Lptraces} holds with $p=\infty$. Note that from the given proof it follows that \eqref{Lptraces} holds for any $p<\infty$. In particular for the smoothing sequence $u^\epsilon=u*\phi^\epsilon \in S_\infty^1(\mathbb{R}^n)\cap C_{loc}^\infty(\mathbb{R}^n)$ all the traces indicated on the left hand side of \eqref{Lptraces} are uniformly bounded in $L_\infty(P_{i_1\cdots i_k})$, and converge to corresponding traces of $u$ in $L_p(P_{i_1\cdots i_k})$ with any $1<p<\infty$. Such limits are also limits in the sense of distributions. Since $L_\infty(P_{i_1\cdots i_k})$ is a dual space of $L_1(P_{i_1\cdots i_k})$, distributional limit of the sequence bonded in $L_\infty(P_{i_1\cdots i_k})$ remains in $L_\infty(P_{i_1\cdots i_k})$. Therefore, \eqref{Lptraces} holds with $p=\infty$.  Theorem is proved.                        $\quad\square$

Corollaries 3 and 4 are direct consequence of the Theorem~\ref{newembedding} due to the fact that if $u\in S_p^k(\mathbb{R}^n$, then all the weak partial derivatives of order $k-1$ are elements of $S_p^1(\mathbb{R}^n)$, and if $u\in S_p^{{\bf k}}(\mathbb{R}^n)$ the indicated partial derivative on the left hand side of \eqref{newembd2} is an element of $S_p^1(\mathbb{R}^n)$. The bounded embedding \eqref{newembddom} is a direct consequence of \eqref{newembd1} and the definition of the space $S_p^1(Q)$. The compact embedding \eqref{newembdcmp} follows from \eqref{newembddom} and Arzela-Ascoli's theorem.


\section{Conclusion}

The concept of Sobolev spaces became a trailblazing idea in many fields of mathematics. The goal of this paper is to gain insight into the embedding of the Sobolev spaces into H\"{o}lder spaces - a very powerful concept that reveals the connection between weak differentiability and integrability (or weak regularity) of the function with its pointwise regularity. It is well-known that the embedding of the Sobolev space of weakly differentiable functions into H\"{o}lder spaces holds if the integrability exponent is higher than the space dimension. Otherwise speaking, one can trade one degree of weak regularity with an integrability exponent higher than the space dimension to upgrade the pointwise regularity to H\"{o}lder continuity. In this paper, the embedding of the Sobolev functions into the H\"{o}lder spaces is expressed in terms of the minimal weak differentiability requirement independent of the integrability exponent. Precisely, the question asked is what is the minimal weak regularity degree of Sobolev functions which upgrades the pointwise regularity to H\"{o}lder continuity independent of the integrability exponent. The paper reveals that the anticipated "largest" subspace of the Sobolev space of weakly differentiable functions embedded into the H\"{o}lder space is the Sobolev space with dominating mixed smoothness. The proof is based on the generalization of the Newton-Leibniz formula to the $n$-dimensional rectangle and inductive application of the Sobolev trace embedding results. The method is applied to prove the embedding of the Sobolev spaces with dominating mixed smoothness into H\"{o}lder spaces.

}
\end{document}